\documentclass{conm-p-l}

\def \cal{\mathcal}
\def \Z{\mathbb Z}
\def \C{\mathbb C}

\def \N{\mathbb N}
\def \wt{{\rm wt}}

\def \Res{{\rm Res}}

\def \Ind {{\rm Ind}}

\def \Hom{{\rm Hom}}
\def \mod{{\rm mod}}

\def \<{\langle} 
\def \>{\rangle}

\def \bconj{\begin{conj}\label}
\def \econj{\end{conj}}
\def \be{\begin{equation}\label}
\def \ee{\end{equation}}
\def \bex{\begin{exa}\label}
\def \eex{\end{exa}}
\def \bl{\begin{lemma}\label}
\def \el{\end{lemma}}
\def \bt{\begin{theorem}\label}
\def \et{\end{theorem}}
\def \bp{\begin{proposition}\label}
\def \ep{\end{proposition}}
\def \br{\begin{remark}\label}
\def \er{\end{remark}}
\def \bc{\begin{corollary}\label}
\def \ec{\end{corollary}}
\def \bd{\begin{definition}\label}
\def \ed{\end{definition}}

\newtheorem{theorem}{Theorem}[section]
\newtheorem{lemma}[theorem]{Lemma}
\newtheorem{proposition}[theorem]{Proposition}
\newtheorem{corollary}[theorem]{Corollary}

\theoremstyle{definition}
\newtheorem{definition}[theorem]{Definition}

\theoremstyle{remark}
\newtheorem{remark}[theorem]{Remark}
\numberwithin{equation}{section}

\makeatletter
\@addtoreset{equation}{section}

\makeatother
\makeatletter
\begin{document}

\title{The regular representations and the $A_{n}(V)$-algebras}
%for vertex operator algebras}

\author{Haisheng Li}
\address{Department of Mathematical Sciences, Rutgers University-Camden, 
Camden, NJ 08102}
%\curraddr{Department of Mathematics and Statistics,
%Case Western Reserve University, Cleveland, Ohio 43403}
\email{hli@crab.rutgers.edu}
\thanks{This research was supported in part 
by NSF grants DMS-9616630 and DMS-9970496.}

%\author{Author Two}
%\address{Mathematical Research Section, School of Mathematical Sciences,
%Australian National University, Canberra ACT 2601, Australia}

%    General info
\date{September 9, 1999.}

%\dedicatory{This paper is dedicated to our authors.}

\begin{abstract}
For a vertex operator algebra $V$, the regular representations
are related to the $A_{n}(V)$-algebras and their bimodules, and 
induced $V$-modules from $A_{n}(V)$-modules are defined and studied
in terms of the regular representations.
\end{abstract}

\maketitle

\section{Introduction}
In [Li2], for a vertex operator algebra $V$ and a nonzero 
complex number $z$,  a weak $V\otimes V$-module
${\cal{D}}_{P(z)}(V)$ was constructed out of the dual space 
$V^{*}$, and certain results of Peter-Weyl type were obtained.
The weak $V\otimes V$-modules ${\cal{D}}_{P(z)}(V)$ were referred
as regular representations.
In [Li3], as a generalization,
weak $V\otimes V$-modules ${\cal{D}}_{P(z)}(V,U)$ were constructed
for any vector space $U$. Furthermore,
Zhu's $A(V)$-theory ([Z1], [FZ]) was related to 
the regular representations in the spirit
of the induced module theory for a Lie group (cf. [Ki]), and 
a notion of an induced $V$-module from an $A(V)$-module 
was formulated in terms of the regular representations. 
The induced $V$-module from an $A(V)$-module $U$
was defined in [Li3] as follows:
First consider linear functions from $V$ to $U$, 
which are lifted from linear functions from $A(V)$ to $U$, 
or simply just linear functions from $A(V)$ to $U$.
Second, it was shown that 
$\Hom(A(V),U)$ is a subspace of
${\cal{D}}_{P(-1)}(V,U)$, and what is more,
$\Hom(A(V),U)$ and 
$\Omega({\cal{D}}_{P(-1)}(V,U))$ $(\subset \Hom(V,U))$
coincide as natural $A(V)\otimes A(V)$-modules. 
Meanwhile, all the (left) $A(V)$-invariant functions from $A(V)$ to $U$
give us the space $\Hom_{A(V)}(A(V),U)$, which is 
canonically isomorphic to $U$ as an $A(V)$-module.
Third, the induced module $\Ind_{A(V)}^{V}U$ was defined
to be the submodule of ${\cal{D}}_{P(-1)}(V,U)$,
generated by $\Hom_{A(V)}(A(V),U)$ $(=U)$ under the action of 
$V\otimes {\C}$.

In [DLM2], as a generalization of Zhu's $A(V)$-theory, a family of associative
algebras $A_{n}(V)$ were constructed and a family of functors 
$\Omega_{n}$ from the category of weak $V$-modules to the
category of $A_{n}(V)$-modules and a family of functors $M_{n}$ 
(with certain properties) from the
category of $A_{n}(V)$-modules to the category of ${\N}$-graded 
weak $V$-modules were constructed.
By definition, $\Omega_{n}(W)$ consists of each $w$ 
such that $v_{m}w=0$ for homogeneous $v\in V$ and for $m\ge \wt v+n$.
(Of course, $\Omega_{n}(W)$ can also be considered as the invariant space
with respect to a certain Lie algebra.)
In the case that $W$ is a lowest weight generalized irreducible $V$-module,
$\Omega_{n}(W)$ is the sum of the first $n$ lowest weight subspaces. 
 
In 1993, Zhu [Z2] gave a general construction of associative algebras
from a vertex operator algebra for a certain purpose. The algebras
$A_{n}(V)$ might be related to those algebras in a certain way.

In this paper, we shall relate $A_{n}(V)$-theory to the 
(generalized) regular representations of $V$ on ${\cal{D}}_{P(-1)}(V,U)$.
When $n\ge 1$, unlike the $n=0$ case [Li3], there are
certain complicated factors. It is proved (Propositions \ref{pan(W)},
\ref{pan'wbimodule}, and Corollary \ref{canvbimodule0})
that as vector spaces, $\Hom (A_{n}(V),U)$
is a subspace of $\Omega_{n}({\cal{D}}_{P(-1)}(V,U))$. However,
as natural $A_{n}(V)\otimes A_{n}(V)$-modules, $\Hom (A_{n}(V),U)$
is not a submodule.
It turns out that the $A_{n}\otimes A_{n}(V)$-module structure
on $\Hom (A_{n}(V),U)$ coincides with a twisted or deformed
$A_{n}\otimes A_{n}(V)$-module structure on 
$\Omega_{n}({\cal{D}}_{P(-1)}(V,U))$
with respect to a certain linear automorphism on 
$\Omega_{n}({\cal{D}}_{P(-1)}(V,U))$
(Theorem \ref{trelation}).
Using this connection we formulate a notion of induced $V$-module from 
an $A_{n}(V)$-module and we show that the induced modules
are lowest weight generalized $V$-modules if the given $A_{n}(V)$-modules 
are irreducible.

An induced module theory from
modules for a vertex operator subalgebra
was established in [DLin]. As mentioned in [Li3], 
the notion of induced module
defined here and the notion of induced module defined in [DLin] are
different in nature.

This paper is organized as follows: In Section 2, we review
the construction of the weak $V\otimes V$-module ${\cal{D}}_{P(z)}(W,U)$.
In Section 3, we relate $A_{n}(V)$-modules
$\Hom(A_{n}(W),U)$ with $\Omega_{n}({\cal{D}}_{P(-1)}(W,U))$, and we define the
induced $V$-module $\Ind_{A_{n}(V)}^{V}U$ for a given $A_{n}(V)$-module $U$.

\section{The weak $V\otimes V$-module ${\cal{D}}_{P(z)}(W,U)$}

In this section we shall recall from [Li3] the construction of the weak 
$V\otimes V$-module ${\cal{D}}_{P(z)}(W,U)$ and there are nothing new.

We use standard definitions and notations as given in [FLM] and [FHL].
A vertex operator algebra is denoted by $(V,Y,{\bf 1},\omega)$,
where ${\bf 1}$ is the vacuum vector and $\omega$ is the Virasoro element,
or simply by $V$.
We also use the notion of weak module as defined in [DLM2]---A weak module
satisfies all the axioms given in [FLM] and [FHL] for the notion 
of a module except that no grading is required.

We typically use letters $x,y, x_{1},x_{2},\dots$ for mutually commuting
formal variables and $z,z_{0},\dots$ for complex numbers.
For a vector space $U$, $U[[x,x^{-1}]]$ is the vector space of all
(doubly infinite) formal series with coefficients in $U$,
$U((x))$ is the space of formal Laurent series in $x$, and
$U((x^{-1}))$ is the space of formal Laurent series in $x^{-1}$.
We emphasize the following standard formal variable convention:
\begin{eqnarray}
& &(x_{1}-x_{2})^{n}=\sum_{i\ge 0}(-1)^{i}\binom{n}{i}x_{1}^{n-i}x_{2}^{i},\\
& &(x-z)^{n}=\sum_{i\ge 0}(-z)^{i}\binom{n}{i}x^{n-i},\\
& &(z-x)^{n}=\sum_{i\ge 0}(-1)^{i}z^{n-i}\binom{n}{i}x^{i}
\end{eqnarray}
for $n\in {\Z},\; z\in {\C}^{\times}$.

For vector spaces $U_{1},U_{2}$, a linear map $f\in \Hom(U_{1},U_{2})$ 
extends canonically to a linear map from $U_{1}[[x,x^{-1}]]$ to 
$U_{2}[[x,x^{-1}]]$. We shall use this canonical extension
without any comments.

Let $V$ be a vertex operator algebra. 
For $v\in V$, we set (cf. [FHL], [HL1])
\begin{eqnarray}
Y^{o}(v,x)=Y(e^{xL(1)}(-x^{-2})^{L(0)}v,x^{-1}).
\end{eqnarray}
For a weak $V$-module $W$, 
because $e^{xL(1)}(-x^{-2})^{L(0)}v\in V[x,x^{-1}]$ and
$Y(u,x^{-1})w\in W((x^{-1}))$ for $u\in V,\; w\in W$,
$Y^{o}(v,x)$ lies in $\Hom (W,W((x^{-1})))$. 
More generally,  for any complex number $z_{0}$,
$Y^{o}(v,x+z_{0})$  lies in $\Hom (W,W((x^{-1})))$, where by definition
\begin{eqnarray}
Y^{o}(v,x+z_{0})w=(Y^{o}(v,y)w)|_{y=x+z_{0}}
\end{eqnarray}
for $w\in W$.
Let $W$ be a weak $V$-module and let $U$ be a vector space, e.g., $U={\C}$.
For $v\in V,\;f\in \Hom(W,U)$, the compositions
$fY^{o}(v,x)$ and $fY^{o}(v,x+z_{0})$ for any complex number $z_{0}$
are elements of $(\Hom (W,U))[[x,x^{-1}]]$. 

Let ${\C}(x)$ be the algebra of rational functions of $x$ (and
${\C}[[x,x^{-1}]]$ be the vector space of all doubly infinite formal
series in $x$ with complex coefficients).
The $\iota$-maps $\iota_{x;0}$ and $\iota_{x;\infty}$ from 
${\C}(x)$ to ${\C}[[x,x^{-1}]]$ are defined as follows:
for any rational function $f(x)$, 
$\iota_{x;0}f(x)$ is the Laurent series expansion of $f(x)$ at $x=0$
and $\iota_{x;\infty}f(x)$ is the Laurent series expansion 
of $f(x)$ at $x=\infty$. These are injective 
${\C}[x,x^{-1}]$-linear maps. In terms of
the formal variable convention, we have
\begin{eqnarray}
& &\iota_{x;0}\left((x-z)^{n}f(x)\right)=(-z+x)^{n}\iota_{x;0}f(x),\\
& &\iota_{x;\infty}\left((x-z)^{n}f(x)\right)=(x-z)^{n}\iota_{x;\infty}f(x)
\end{eqnarray}
for $n\in {\Z},\; z\in {\C}^{\times},\; f(x)\in {\C}(x)$.

\bd{dDWU}
Let $W$ be a weak $V$-module, $U$ a vector space and 
$z$ a nonzero complex number.
Define ${\cal{D}}_{P(z)}(W,U)$ to be the subspace of $\Hom(W,U)$,
consisting of each $f$ such that
for $v\in V$, there exist $k,l\in {\N}$ 
such that
\begin{eqnarray}\label{eDWUcharc}
(x-z)^{k}x^{l}\<u^{*}, fY^{o}(v,x)w\>\in {\C}[x]
\end{eqnarray}
for all $u^{*}\in U^{*},\; w\in W$, or what is equivalent,
for all $u^{*}\in U^{*},\; w\in W$, 
the formal series
$$\<u^{*}, fY^{o}(v,x)w\>,$$
an element of ${\C}((x^{-1}))$,
absolutely converges in the domain $|x|>|z|$ to a rational function
of the form $x^{-l}(x-z)^{-k}g(x)$ for $g(x)\in {\C}[x]$.
\ed

The following are equivalent definitions of ${\cal{D}}_{P(z)}(W,U)$
in terms of formal series:

\bl{ldefdwuequiv} 
Let $f\in \Hom(W,U)$. Then the following statements are equivalent:

(a) $f\in {\cal{D}}_{P(z)}(W,U)$.

(b) For $v\in V$, there exist $k,l\in {\N}$ such that
\begin{eqnarray}\label{edefDWU}
(x-z)^{k}x^{l}fY^{o}(v,x)\in (\Hom(W,U))[[x]].
\end{eqnarray}

(c)  For $v\in V$, there exist $k,l\in {\N}$ such that for each $w\in W$,
\begin{eqnarray}\label{edefDWU2}
(x-z)^{k}x^{l}fY^{o}(v,x)w\in U[x].
\end{eqnarray}
\el

Let $v\in V,\; f\in {\cal{D}}_{P(z)}(W,U)$ and let $k,l\in {\N}$
be such that (\ref{edefDWU2}) holds.
Then by changing variable we get
\begin{eqnarray}\label{edefDWU3}
x^{k}(x+z)^{l}fY^{o}(v,x+z)w\in U[x]
\end{eqnarray}
for $w\in W$.

\bd{dLRactions}
Let $W, U$ and $z$ be given as before. 
For 
$$v\in V,\;\;\;\; f\in {\cal{D}}_{P(z)}(W,U),$$
we define
two elements $Y^{L}_{P(z)}(v,x)f$ and $Y_{P(z)}^{R}(v,x)f$
of $(\Hom (W,U))[[x,x^{-1}]]$ by
\begin{eqnarray}
(Y_{P(z)}^{L}(v,x)f)(w)
&=&(z+x)^{-l}\left((x+z)^{l}f(Y^{o}(v,x+z)w)\right)\\
(Y_{P(z)}^{R}(v,x)f)(w)&=&(-z+x)^{-k}\left((x-z)^{k}f(Y^{o}(v,x)w)\right)
\end{eqnarray}
for $w\in W$, where $k,l$ are any pair of (possibly negative) integers 
such that (\ref{edefDWU}) holds.
\ed

First, in view of (\ref{edefDWU2}) and (\ref{edefDWU3}), 
both $(z+x)^{-l}\left((x+z)^{l}f(Y^{o}(v,x+z)w)\right)$
and $(-z+x)^{-k}\left((x-z)^{k}f(Y^{o}(v,x)w)\right)$ lie
in $U((x))$, so that 
$Y^{L}_{P(z)}(v,x)f$ and $Y^{R}_{P(z)}(v,x)f$ make sense.
However, we are not allowed to
remove the left-right brackets to cancel $(x-z)^{k}$ or $(x+z)^{l}$ 
because of the nonexistence of $(z+x)^{-l}f(Y^{o}(v,x+z)w)$ and 
$(-z+x)^{-k}f(Y^{o}(v,x)w)$.
Second, they are also well defined, i.e., they are independent
of the choice of the pair of integers $k,l$. Indeed, 
if $k',l'$ are another pair of integers such that (\ref{edefDWU}) holds,
say for example, $k\ge k'$, then 
\begin{eqnarray}
& &(-z+x)^{-k}\left((x-z)^{k}fY^{o}(v,x)w\right)\\
&=&(-z+x)^{-k}\left((x-z)^{k-k'}(x-z)^{k'}fY^{o}(v,x)w\right)\nonumber\\
&=&(-z+x)^{-k}(x-z)^{k-k'}\left((x-z)^{k'}fY^{o}(v,x)w\right)\nonumber\\
&=&(-z+x)^{-k'}\left((x-z)^{k'}fY^{o}(v,x)w\right).\nonumber
\end{eqnarray}

{}From definition we immediately have:

\bl{lDWUconn}
For $v\in V,\; f\in {\cal{D}}_{P(z)}(W,U)$,
\begin{eqnarray}
& &(z+x)^{l}Y_{P(z)}^{L}(v,x)f=(x+z)^{l}fY^{o}(v,x+z),\\
& &(-z+x)^{k}Y_{P(z)}^{R}(v,x)f=(x-z)^{k}fY^{o}(v,x),
\end{eqnarray}
where $k,l$ are any pair of (maybe negative) integers 
such that (\ref{edefDWU}) holds.
\el

{}From the definition, 
$\<u^{*},fY^{o}(v,x)w\>$ lies in the range of $\iota_{x;\infty}$
for $u^{*}\in U^{*},\;f\in {\cal{D}}_{P(z)}(W,U),\; v\in V,\; w\in W$.
Then $\iota_{x;\infty}^{-1}\<u^{*}, fY^{o}(v,x)w\>$ is a well 
defined element of ${\C}(x)$.
In terms of rational functions and the $\iota$-maps we immediately have:

\bl{lDWUconn2}
For $v\in V,\; f\in {\cal{D}}_{P(z)}(W,U),\; u^{*}\in U^{*},\; w\in W$,
\begin{eqnarray}
& &\<u^{*},(Y_{P(z)}^{L}(v,x)f)(w)\>
=\iota_{x;0}\iota_{x;\infty}^{-1}\<u^{*},fY^{o}(v,x+z)w\>,\\
& &\<u^{*},(Y_{P(z)}^{R}(v,x)f)(w)\>
=\iota_{x;0}\iota_{x;\infty}^{-1}\<u^{*},fY^{o}(v,x)w\>.
\end{eqnarray}
\el

\bt{tDWU}
Let $W$ be a weak $V$-module, $U$ a vector space and $z$ a nonzero 
complex number. Then the pairs $({\cal{D}}_{P(z)}(W,U), Y_{P(z)}^{L})$ and
$({\cal{D}}_{P(z)}(W,U), Y_{P(z)}^{R})$ carry 
the structure of a weak $V$-module and the actions
$Y_{P(z)}^{L}$ and $Y_{P(z)}^{R}$ of $V$ on ${\cal{D}}_{P(z)}(W,U)$
commute. Furthermore, set 
\begin{eqnarray}
Y_{P(z)}=Y_{P(z)}^{L}\otimes Y_{P(z)}^{R}.
\end{eqnarray}
Then the pair $({\cal{D}}_{P(z)}(W,U), Y_{P(z)})$ 
carries the structure of a weak $V\otimes V$-module. 
\et

The following 
relation among $fY^{o}(v,x), Y^{L}(v,x)f$ and $Y^{R}(v,x)f$ holds [Li3]:

\bp{pthreeconn}
Let $v\in V,\; f\in {\cal{D}}_{P(z)}(W,U)$. Then
\begin{eqnarray}\label{ethreeconn}
& &x_{0}^{-1}\delta\left(\frac{x-z}{x_{0}}\right)fY^{o}(v,x)
-x_{0}^{-1}\delta\left(\frac{z-x}{-x_{0}}\right)Y_{P(z)}^{R}(v,x)f\\
&=&z^{-1}\delta\left(\frac{x-x_{0}}{z}\right)Y_{P(z)}^{L}(v,x_{0})f.\nonumber
\end{eqnarray}
\ep

\section{The associative algebras $A_{n}(V)$ and induced modules
 $\Ind_{A_{n}(V)}^{V}U$}

In this section, the nonzero complex number $z$ in the notion of weak 
$V\otimes V$-module ${\cal{D}}_{P(z)}(W,U)$ will be specified as
$-1$. We shall simply use $Y^{L}$ and $Y^{R}$ for
$Y^{L}_{P(-1)}$ and $Y^{R}_{P(-1)}$.
Throughout this section, $n$ will represent a nonnegative integer.

We shall need the following notions.
A {\em generalized} $V$-module [HL1] is a weak $V$-module on which
$L(0)$ semisimply acts. 
 Then for a generalized $V$-module $W$ we have 
the $L(0)$-eigenspace decomposition: $W=\coprod_{h\in {\C}}W_{(h)}$.
Thus, a generalized $V$-module 
satisfies all the axioms defining the notion of a $V$-module
([FLM], [FHL]) except the two grading restrictions on 
homogeneous subspaces.
If a generalized $V$-module furthermore
satisfies the lower truncation condition (one of the two grading 
restrictions), it is called a {\em lower truncated} generalized module [H1].

A {\em lowest weight} generalized $V$-module is a generalized $V$-module 
such that $W=\coprod_{n\in \N}W_{(h+n)}$ for some $h\in {\C}$ and
$W_{(h)}$ generates $W$ as a weak $V$-module.
Furthermore, if $W\ne 0$, we call the unique $h$ 
{\em the lowest weight} of $W$. 

Now we recall the construction of $A_{n}(V)$ algebra and some basic 
results from [DLM2].

\bd{donw} Let $V$ be a vertex operator algebra and let $n\in {\N}$.
Define a subspace $O_{n}(V)$ of $V$, linearly spanned by elements
\begin{eqnarray}
& &(L(-1)+L(0))v,\\
& &{\rm Res}_{x}\frac{(1+x)^{\wt u+n}}{x^{2n+2}}Y(u,x)v
\end{eqnarray}
for homogeneous $u,v\in V$.
Define
\begin{eqnarray}\label{eu*nv}
u*_{n}v&=&\sum_{m=0}^{n}(-1)^{m}\binom{m+n}{n}{\rm Res}_{x}
\frac{(1+x)^{\wt u+n}}{x^{n+m+1}}Y(u,x)v\\
& &\left(=\sum_{m=0}^{n}\binom{-n-1}{m}{\rm Res}_{x}
\frac{(1+x)^{\wt u+n}}{x^{n+m+1}}Y(u,x)v\right).\nonumber
\end{eqnarray}
\ed

We have:

\bl{lrightactionanv} {\rm [DLM2]}
Let $u,v\in V$ be homogeneous. Then
 \begin{eqnarray*}
u*_{n}v-\sum_{m=0}^{n}\binom{-n-1}{m}(-1)^{n-m}{\rm Res}_{x}
x^{-n-m-1}(1+x)^{\wt v+m-1}Y(v,x)u\in O_{n}(V).
\end{eqnarray*}
\el

Set $A_{n}(V)=V/O_{n}(V)$.

\bp{pdlm2} {\rm [DLM2]} Let $(V,Y,{\bf 1},\omega)$ be a 
vertex operator algebra. Then

(a) $O_{n}(V)$ is a two-sided ideal of the nonassociative algebra $(V,*_{n})$
and the quotient algebra
$A_{n}(V)$ is an associative algebra with ${\bf 1}+O_{n}(V)$ 
as its identity element, with $\omega+O_{n}(V)$ being central  and with an 
involution (anti-automorphism)
\begin{eqnarray}\label{etheta}
\theta: v+O_{n}(V)\mapsto e^{L(1)}(-1)^{L(0)}v+O_{n}(V).
\end{eqnarray}

(b) For each $n\ge 0$, the identity map of $V$ gives rise to an algebra 
homomorphism $\psi_{n}$ from $A_{n+1}(V)$ onto $A_{n}(V)$.
\ep

For any weak $V$-module $W$ and any $n\in {\N}$, we define [DLM2]
\begin{eqnarray}\label{eomeganw}
\mbox{}\;\;\;\;\;\;\Omega_{n}(W)=\{ w\in W\;|\; v_{\wt v+m}w=0\;\;\;
\mbox{ for homogeneous }v\in V,\; m\ge n\}.
\end{eqnarray}
For $w\in \Omega_{n}(W)$ and for homogeneous $v\in V$, we have 
$x^{\wt v+n}Y(v,x)w\in W[[x]]$.
When $n=0$, we have $\Omega(W)=\Omega_{0}(W)$. Clearly,
\begin{eqnarray}
\Omega_{0}(W)\subset \Omega_{1}(W)\subset \cdots.
\end{eqnarray}
We similarly define $\Omega_{-1}(W), \Omega_{-2}(W),\dots$.
Since $\wt {\bf 1}=0$ and ${\bf 1}_{r}w=\delta_{r,-1}w$, 
we have $\Omega_{-n}(W)=0$ for $n\ge 1$.

The following result was proved in [DLM2]:

\bp{pomeganw}
Let $W$ be a weak $V$-module and let $n\ge 0$. Then
$\Omega_{n}(W)$ is an $A_{n}(V)$-module where $v+O_{n}(V)$
acts as $v_{\wt v-1}$ for homogeneous $v\in V$.
\ep

Let $W_{1}, W_{2}$ be weak $V$-modules and let $\psi$ be a 
$V$-homomorphism from $W_{1}$ to $W_{2}$. It is clear that
$\psi(\Omega_{n}(W_{1}))\subset \Omega_{n}(W_{2})$ and
the restriction $\Omega_{n}(\psi):=\psi|_{\Omega_{n}(W_{1})}$ is an 
$A_{n}(V)$-homomorphism. It is routine to check that
we have obtained a functor $\Omega_{n}$ from 
the category of weak $V$-modules to the category of $A_{n}(V)$-modules.

\bl{lomegan}
Let $W$ be a weak $V$-module and set
\begin{eqnarray}
{\cal{S}}(W)=\cup_{n\ge 0}\Omega_{n}(W).
\end{eqnarray}
Let $u\in V$ be homogeneous and let $r\in {\Z}$.
Then 
\begin{eqnarray}
u_{r}\Omega_{n}(W)\subset \Omega_{n}(W)
\end{eqnarray}
if $r\ge \wt u-1$, and 
\begin{eqnarray}
u_{r}\Omega_{n}(W)\subset \Omega_{n+\wt u-r-1}(W)
\end{eqnarray}
if $r<\wt u-1$.
In particular, ${\cal{S}}(W)$ is a sub-weak-module of $W$.
Furthermore,
\begin{eqnarray}
\Omega_{n}({\cal S}(W))=\Omega_{n}(W).
\end{eqnarray}
\el

\begin{proof} Let $w\in \Omega_{n}(W)$, let $v\in V$ be homogeneous and 
let $m\in {\Z}$. By Borcherds commutator formula,
\begin{eqnarray}
& &v_{\wt v+m}u_{r}w\\
&=&u_{r}v_{\wt v+m}w
+\sum_{i\ge 0}\binom{\wt v+m}{i}(v_{i}u)_{\wt v+m+r-i}w
\nonumber\\
&=&u_{r}v_{\wt v+m}w+\sum_{i\ge 0}\binom{\wt v+m}{i}
(v_{i}u)_{\wt (v_{i}u)+m+r-\wt u+1}w.\nonumber
\end{eqnarray}
Then the first part follows immediately.
Since ${\cal S}(W)$ is a submodule of $W$, we have 
$\Omega_{n}({\cal S}(W))\subset \Omega_{n}(W)$. 
It is easy to see that
 $\Omega_{n}(W)\subset \Omega_{n}({\cal S}(W))$. This completes the proof.
\end{proof}

We shall need the result that $L(1)$ is locally nilpotent on 
${\cal S}(W)$ for any weak $V$-module $W$. To prove this result,
we recall from [Li3] the following result, which is a reformulation
of a result in [DLM2] (Remark 3.3):

\bl{lassoc}
Let $W$ be a weak $V$-module, $w\in W$. Let $u,v\in V$ and
let $k\in {\Z}$ be such that
\begin{eqnarray}
x^{k}Y(u,x)w\in W[[x]],
\end{eqnarray}
or equivalently,
\begin{eqnarray}
u_{k+m}w=0\;\;\mbox{ for }m\ge 0.
\end{eqnarray}
Then for $p,q\in {\Z}$,
\begin{eqnarray}
u_{p}v_{q}w
=\sum_{i=0}^{s}\sum_{j\ge 0}\binom{p-k}{i}\binom{k}{j}
(u_{p-k-i+j}v)_{q+k+i-j}w.
\end{eqnarray}
where $s$ is any nonnegative integer such that $x^{s+1+q}Y(v,x)w\in W[[x]]$. 
\el

As an immediate consequence we have ([DM] and [Li1]):

\bc{cdmli}
Let $W$ be a weak $V$-module and let $w\in W$. Set
\begin{eqnarray}
\<w\>=\mbox{{\rm linear span }}\{ v_{m}w\;|\; v\in V,\; m\in {\Z}\}.
\end{eqnarray}
Then $\<w\>$ is the sub-weak-module of $W$, generated by $w$.
\ec

\bl{lnilpotentL(1)}
Let $W$ be a weak $V$-module. Then for any $r$
homogeneous vectors $v^{1},\dots,v^{r}\in V$,
\begin{eqnarray}
v^{1}_{m_{1}}\cdots v^{r}_{m_{r}}\Omega_{n}(W)=0
\end{eqnarray}
for  $m_{i}\in {\Z}$ with
$$m_{1}+\cdots + m_{r}\ge \wt v^{1}+\cdots +\wt v^{r}-r +n.$$ 
In particular, for homogeneous $v\in V$ and for $m\ge \wt v$, 
\begin{eqnarray}
(v_{m})^{n}\Omega_{n}(W)=0.
\end{eqnarray}
\el

\begin{proof}
 We shall prove the first part by induction on $r$. From the definition 
of $\Omega_{n}(W)$, 
the lemma is true for $r=1$. Assume it is true for any $r$ homogeneous 
vectors in $V$. Now let $v^{1},\dots,v^{r}, v^{r+1}\in V$ be homogeneous
and let $m_{i}\in {\Z}$ with
\begin{eqnarray}
m_{1}+\cdots + m_{r}+m_{r+1}\ge \wt v^{1}+\cdots +\wt v^{r+1}-(r+1) +n.
\end{eqnarray}
Set 
$$u=v^{r},\;\; v=v^{r+1},\;\; p=m_{r},\;\; q=m_{r+1}.$$
Since $w\in \Omega_{n}(W)$, in Lemma \ref{lassoc}, we may take
$k=\wt u+n=\wt v^{r}+n$. Let $s$ be any nonnegative integer
such that $x^{s+1+q}Y(v,x)w\in W[[x]]$. By Lemma \ref{lassoc},
we have
\begin{eqnarray}
& &u_{p}v_{q}w\\
&=&\sum_{i=0}^{s}\sum_{j\ge 0}\binom{p-\wt u-n}{i}\binom{\wt u+n}{j}
(u_{p-\wt u-n-i+j}v)_{q+\wt u+n+i-j}w.\nonumber
\end{eqnarray}
Notice that
\begin{eqnarray}
\wt (u_{p-\wt u-n-i+j}v)&=&\wt u+\wt v+\wt u+n+i-j-1-p\\
&=&2\wt v^{r}+\wt v^{r+1}+n+i-j-1-p.\nonumber
\end{eqnarray}
Thus
\begin{eqnarray}
& &m_{1}+\cdots +m_{r-1}+(q+\wt u+n+i-j)\nonumber\\
&\ge& \wt v^{1}+\cdots +\wt v^{r+1}-(r+1)+n+(q+\wt u+n+i-j)-m_{r}-m_{r+1}
\nonumber\\
&=&\wt v^{1}+\cdots +\wt v^{r-1}+\wt (u_{p-\wt u-n-i+j}v)-r+n.\nonumber
\end{eqnarray}
Then it follows from the inductive hypothesis that
\begin{eqnarray}
& &v^{1}_{m_{1}}\cdots v^{k+1}_{m_{k+1}}w\\
&=&\sum_{i=0}^{s}\sum_{j\ge 0}\binom{p-k}{i}\binom{k}{j}
v^{1}_{m_{1}}\cdots v^{k-1}_{m_{k-1}}(u_{p-k-i+j}v)_{q+k+i-j}w\nonumber\\
&=&0.\nonumber
\end{eqnarray}
This finishes the induction and concludes the proof.
\end{proof}

In view of Lemma \ref{lnilpotentL(1)}, 
noticing that $L(1)=\omega_{2}$ and $\wt \omega=2$, 
we immediately have:

\bc{cnilpotentL(1)}
Let $W$ be a weak $V$-module, let $v\in V$ be homogeneous and let $m\ge \wt v$.
Then $v_{m}$ is locally nilpotent on ${\cal{S}}(W)$.
In particular, $L(1)$ is locally nilpotent on ${\cal{S}}(W)$.
\ec

Let $W$ be a weak $V$-module. We define
$O'_{n}(W)$ to be the subspace of 
$W$, linearly spanned by elements of the form:
\begin{eqnarray}
v\circ_{n}w:=\Res_{x}x^{-2n-2}(1+x)^{\wt v+n}Y(v,x)w
\end{eqnarray}
for $w\in W$ and for homogeneous $v\in V$. 
The proof of Lemma 2.1.2 of [Z1] with minor necessary changes directly gives:

\bl{lzhuO'(W)}
Let $W$ be a weak $V$-module, let $w\in W$, and let $v\in V$ be homogeneous.
Then
\begin{eqnarray}
\Res_{x}x^{-2n-2-r}(1+x)^{\wt v+n+s}Y(v,x)w\in O'_{n}(W)
\end{eqnarray}
for $r\ge s\ge 0$.
\el

In the following there will be several module structures 
on a certain vector space.
For this reason,  we shall use $\Omega_{n}(W,Y_{W})$ 
including the vertex operator map $Y_{W}$ in the notation 
for $\Omega_{n}(W)$. Since $\Hom(-,U)$
is a contravariant functor for the category of vector spaces,
for any vector spaces $A, B$ and any surjective linear map $g\in \Hom(A,B)$, 
we have 
an injective linear map $\Hom(g,U)$ from $\Hom (B,U)$ into $\Hom(A,U)$.
In particular, if $B$ is a quotient space of $A$, we may naturally
identify $\Hom (B,U)$ as a subspace of $\Hom(A,U)$.

\bp{pan(W)}
Let $W$ be a weak $V$-module and let $U$ be a vector space.
Set $A'_{n}(W)=W/O'_{n}(W)$. Then
\begin{eqnarray}
\mbox{}\;\;\;\;\;\;\;\Hom (A'_{n}(W),U)
=\Omega_{n}({\cal{D}}_{P(-1)}(W,U),Y^{L})
\cap \Omega_{n}({\cal{D}}_{P(-1)}(W,U),Y^{R}).
\end{eqnarray}
Furthermore,  elements $\alpha$ of
$\Hom(A'_{n}(W),U)$, a natural subspace of $\Hom (W,U)$, are 
characterized by the following property:
\begin{eqnarray}\label{echarcter}
x^{\wt v+n}(x+1)^{\wt v+n}\alpha Y^{o}(v,x)\in (\Hom (W,U))[[x]]
\end{eqnarray}
for homogeneous $v\in V$.
\ep

\begin{proof}
 Let $T$ be the set defined by the property (\ref{echarcter}).
We shall prove
\begin{eqnarray}
& &\Hom (A'_{n}(W),U)\subset T\subset \Omega_{n}({\cal{D}}_{P(-1)}(W,U),Y^{L})
\cap \Omega_{n}({\cal{D}}_{P(-1)}(W,U),Y^{R})\subset \nonumber\\
& &\subset T\subset \Hom (A'_{n}(W),U).\nonumber
\end{eqnarray}
Let $\alpha\in \Hom (A'_{n}(W),U)$ and let $v\in V$ be homogeneous. 
Then for any $m\ge 0$,
\begin{eqnarray}
& &\Res_{x}x^{\wt v+n+m}(x+1)^{\wt v+n}\alpha Y^{o}(v,x)w\\
&=&\Res_{x}x^{\wt v+n+m}(x+1)^{\wt v+n}
\alpha Y(e^{xL(1)}(-x^{-2})^{L(0)}v,x^{-1})w\nonumber\\
&=&(-1)^{\wt v}\Res_{x}x^{2n+m}(1+x^{-1})^{\wt v+n}
\alpha Y(e^{xL(1)}v,x^{-1})w\nonumber\\
&=&(-1)^{\wt v}\Res_{x}x^{-2n-m-2}(1+x)^{\wt v+n}
\alpha Y(e^{x^{-1}L(1)}v,x)w\nonumber\\
&=&0\nonumber
\end{eqnarray}
because (Lemma \ref{lzhuO'(W)})
\begin{eqnarray}
& &\Res_{x}x^{-2n-m-2}(1+x)^{\wt v+n}Y(e^{x^{-1}L(1)}v,x)w\\
&=&\sum_{i\ge 0}\frac{1}{i!}
\Res_{x}x^{-2n-m-2-i}(1+x)^{\wt (L(1)^{i}v)+n+i}Y(L(1)^{i}v,x)w\nonumber\\
&\in& O'_{n}(W).\nonumber
\end{eqnarray}
This proves (\ref{echarcter}). Since $Y^{o}(v,x)w\in W((x^{-1}))$ for $w\in W$,
(\ref{echarcter}) implies
\begin{eqnarray}
x^{\wt v+n}(x+1)^{\wt v+n}\alpha Y^{o}(v,x)w\in U[x].
\end{eqnarray}
By changing variable we get
\begin{eqnarray}
(x-1)^{\wt v+n}x^{\wt v+n}\alpha Y^{o}(v,x-1)\in (\Hom(W,U))[[x]].
\end{eqnarray}
By Lemma \ref{ldefdwuequiv},
$\alpha\in {\cal{D}}_{P(-1)}(W,U)$ and by Lemma \ref{lDWUconn}
\begin{eqnarray}
& &x^{\wt v+n}(1+x)^{\wt v+n}Y^{R}(v,x)\alpha
=x^{\wt v+n}(x+1)^{\wt v+n}\alpha Y^{o}(v,x),\\
& &(-1+x)^{\wt v+n}x^{\wt v+n}Y^{L}(v,x)\alpha
=(x-1)^{\wt v+n}x^{\wt v+n}\alpha Y^{o}(v,x-1).
\end{eqnarray}
Consequently,
\begin{eqnarray}\label{ealphainomega2}
x^{\wt v+n}Y^{R}(v,x)\alpha,\;\; x^{\wt v+n}Y^{L}(v,x)\alpha
\in (\Hom (W,U))[[x]].
\end{eqnarray}
That is, 
$\alpha \in \Omega_{n}({\cal{D}}_{P(-1)}(W,U),Y^{L})
\cap \Omega_{n}({\cal{D}}_{P(-1)}(W,U),Y^{R}).$

Conversely, let $\alpha \in \Omega_{n}({\cal{D}}_{P(-1)}(W,U),Y^{L})
\cap \Omega_{n}({\cal{D}}_{P(-1)}(W,U),Y^{R}).$
Then (\ref{ealphainomega2}) holds
for each homogeneous $v\in V$. 
Recall (\ref{ethreeconn}) with $f=\alpha,\; z=-1$:
\begin{eqnarray}\label{sjacnew}
& &x_{0}^{-1}\delta\left(\frac{x+1}{x_{0}}\right)
\alpha Y^{o}(v,x)-x_{0}^{-1}\delta\left(\frac{1+x}{x_{0}}\right)
Y_{P(z)}^{R}(v,x)\alpha\\
&=&-\delta(-x+x_{0})Y_{P(z)}^{L}(v,x_{0})\alpha.\nonumber
\end{eqnarray}
Applying $\Res_{x_{0}}x^{\wt v+n}x_{0}^{\wt v+n}$ to (\ref{sjacnew}),
then using (\ref{ealphainomega2}) and 
the fundamental properties of delta function we get
\begin{eqnarray}
& &x^{\wt v+n}(x+1)^{\wt v+n}\alpha Y^{o}(v,x)\\
&=&x^{\wt v+n}(1+x)^{\wt v+n}Y^{R}(v,x)\alpha\nonumber\\
& &-\Res_{x_{0}}x^{\wt v+n}x_{0}^{\wt v+n}
\delta(-x+x_{0})Y^{L}(v,x_{0})\alpha\nonumber\\
&=&x^{\wt v+n}(1+x)^{\wt v+n}Y^{R}(v,x)\alpha\nonumber\\
&\in& (\Hom (W,U))[[x]].\nonumber
\end{eqnarray}
Furthermore, for any $w\in W$, 
\begin{eqnarray*}
& &\Res_{x}x^{-2n-2}(1+x)^{\wt v+n}\alpha Y(v,x)w\\
&=&\Res_{x}x^{-2n-2}(1+x)^{\wt v+n}
\alpha Y^{o}(e^{xL(1)}(-x^{-2})^{L(0)}v,x^{-1})w\nonumber\\
&=&\Res_{x}x^{2n}(1+x^{-1})^{\wt v+n}
\alpha Y^{o}(e^{x^{-1}L(1)}(-x^{2})^{L(0)}v,x)w\nonumber\\
&=&\Res_{x}(-1)^{\wt v}x^{\wt v+n}(x+1)^{\wt v+n}
\alpha Y^{o}(e^{x^{-1}L(1)}v,x)w\nonumber\\
&=&\sum_{i\ge 0}(-1)^{\wt v}\frac{1}{i!}\Res_{x}x^{\wt v-i+n}(x+1)^{\wt v+n}
\alpha Y^{o}(L(1)^{i}v,x)w\nonumber\\
&=&\sum_{i\ge 0}(-1)^{\wt v}\frac{1}{i!}
\Res_{x}x^{\wt (L(1)^{i}v)+n}(x+1)^{\wt (L(1)^{i}v)+n+i}
\alpha Y^{o}(L(1)^{i}v,x)w\nonumber\\
&=&0.\nonumber
\end{eqnarray*}
Thus $\alpha(O'_{n}(W))=0$, hence $\alpha\in \Hom(A'_{n}(W),U)$.
This completes the proof.
\end{proof}

It follows from Theorem \ref{tDWU} and 
Proposition \ref{pomeganw} that $\Omega_{n}({\cal{D}}_{P(-1)}(W,U),Y^{R})$
is a natural $A_{n}(V)$-module. Since $Y^{L}$ and $Y^{R}$ commute,
$\Omega_{n}({\cal{D}}_{P(-1)}(W,U),Y^{R})$ is also a weak $V$-module
under the vertex operator map $Y^{L}$.
Then it follows from Proposition \ref{pomeganw} again that
$$\Omega_{n}\left(\Omega_{n}({\cal{D}}_{P(-1)}(W,U),Y^{R}),Y^{L}\right)$$
is an $A_{n}(V)\otimes A_{n}(V)$-module.
Clearly,
\begin{eqnarray}
& &\Omega_{n}\left(\Omega_{n}({\cal{D}}_{P(-1)}(W,U),Y^{R}),Y^{L}\right)\\
&=&\Omega_{n}({\cal{D}}_{P(-1)}(W,U),Y^{L})\cap 
\Omega_{n}({\cal{D}}_{P(-1)}(W,U),Y^{R}).\nonumber
\end{eqnarray}
Thus, $\Omega_{n}({\cal{D}}_{P(-1)}(W,U),Y^{L})\cap 
\Omega_{n}({\cal{D}}_{P(-1)}(W,U),Y^{R})$ is an 
$A_{n}(V)\otimes A_{n}(V)$-module.
For convenience, we refer to this $A_{n}(V)\otimes A_{n}(V)$-module
structure as the {\em canonical module structure}.
{}From definition, we have
\begin{eqnarray}\label{eanv=intersection}
\mbox{}\;\;\;\;\;\;\;\Omega_{n}({\cal{D}}_{P(-1)}(W,U))\subseteq 
\Omega_{n}({\cal{D}}_{P(-1)}(W,U),Y^{L})\cap 
\Omega_{n}({\cal{D}}_{P(-1)}(W,U),Y^{R}).
\end{eqnarray}
(The equality of (\ref{eanv=intersection}) holds when $n=0$, 
but the equality does not hold for $n\ge 1$.) 
It is easy to see that $\Omega_{n}({\cal{D}}_{P(-1)}(W,U))$ is an
$A_{n}(V)\otimes A_{n}(V)$-submodule.

Motivated by [Li3] for $n=0$, we should 
identify $\Hom (A_{n}'(W),U)$ with 
$$\Omega_{n}({\cal{D}}_{P(-1)}(W,U),Y^{L})\cap 
\Omega_{n}({\cal{D}}_{P(-1)}(W,U),Y^{R})$$
as natural $A_{n}(V)\otimes A_{n}(V)$-modules.
We shall prove that $A_{n}'(W)$ just like $A_{n}(V)$ has 
a natural 
$A_{n}(V)\otimes A_{n}(V)$-module structure and so does $\Hom (A_{n}'(W),U)$.
It turns out that the $A_{n}(V)\otimes A_{n}(V)$-module 
$\Hom (A_{n}'(W),U)$ is naturally isomorphic to
$$\Omega_{n}({\cal{D}}_{P(-1)}(W,U),Y^{L})\cap 
\Omega_{n}({\cal{D}}_{P(-1)}(W,U),Y^{R})$$
with a deformed $A_{n}(V)\otimes A_{n}(V)$-module structure.

To achieve our goal we shall need the following result 
(cf. [Li2], Remark 2.10):

\bp{pdeformz}
Let $(E,Y_{E})$ be a weak $V$-module on which $L(1)$ is locally nilpotent,
and let $z_{0}$ be any complex number.
For $v\in V$, we define
\begin{eqnarray}
Y_{E}^{[z_{0}]}(v,x)=Y_{E}(e^{-z_{0}(1+z_{0}x)L(1)}(1+z_{0}x)^{-2L(0)}v,
x/(1+z_{0}x)).
\end{eqnarray}
Then the pair $(E,Y_{E}^{[z_{0}]})$ carries the structure of a weak $V$-module
and $e^{-z_{0}L(1)}$ is a $V$-isomorphism from $(E,Y_{E})$ to 
$(E,Y_{E}^{[z_{0}]})$. 
Furthermore, for homogeneous $v\in V$ and for $m\in {\Z}$, we have
\begin{eqnarray}\label{edeformzrelation}
& &\Res_{x}x^{m}Y_{E}^{[z_{0}]}(v,x)\\
&=&\Res_{x}x^{m}(1-z_{0}x)^{2\wt v-m-2}
Y_{E}(e^{-z_{0}(1-z_{0}x)^{-1}L(1)}v,x).\nonumber
\end{eqnarray}
In particular, 
\begin{eqnarray}\label{edefylw0}
& &\Res_{x}x^{\wt v-1}Y^{[z_{0}]}(v,x)w\\
&=&\Res_{x}x^{\wt v-1}(1-z_{0}x)^{\wt v-1}
Y\left(e^{-z_{0}(1-z_{0}x)^{-1}L(1)}v,x\right).\nonumber
\end{eqnarray}
\ep

\begin{proof}
 Recall the conjugation formula (5.2.38) of [FHL]:
\begin{eqnarray}\label{efhlconjugation}
& &e^{-x_{1}L(1)}Y(v,x)e^{x_{1}L(1)}\\
&=&Y(e^{-x_{1}(1+x_{1}x)L(1)}(1+x_{1}x)^{-2L(0)}v,x/(1+x_{1}x)).\nonumber
\end{eqnarray}
Because $L(1)$ is locally nilpotent on $E$, we may set $x_{1}=z_{0}$,
so that we have
\begin{eqnarray}
& &e^{-z_{0}L(1)}Y_{E}(v,x)e^{z_{0}L(1)}\\
&=&Y_{E}(e^{-z_{0}(1+z_{0}x)L(1)}(1+z_{0}x)^{-2L(0)}v,x/(1+z_{0}x))\nonumber\\
&=&Y_{E}^{[z_{0}]}(v,x).\nonumber
\end{eqnarray}
Then the first part of the proposition follows immediately.

By changing variable $x=y/(1-z_{0}y)$ we get
\begin{eqnarray*}
& &\Res_{x}x^{m}Y_{E}^{[z_{0}]}(v,x)\\
&=&\Res_{x}x^{m}Y_{E}(e^{-z_{0}(1+z_{0}x)L(1)}(1+z_{0}x)^{-2L(0)}v,
x/(1+z_{0}x))\nonumber\\
&=&\Res_{y}y^{m}(1-z_{0}y)^{-m-2}
Y_{E}(e^{-z_{0}(1+z_{0}y)^{-1}L(1)}(1-z_{0}y)^{2L(0)}v,y)\nonumber\\
&=&\Res_{y}y^{m}(1-z_{0}y)^{2\wt v-m-2}
Y_{E}(e^{-z_{0}(1-z_{0}y)^{-1}L(1)}v,y).\nonumber
\end{eqnarray*}
This completes the proof. 
\end{proof}

By definition we have
\begin{eqnarray}
& &(Y^{[z_{0}]})^{[-z_{0}]}(v,x)\\
&=&Y^{[z_{0}]}(e^{z_{0}(1-z_{0}x)L(1)}(1-z_{0}x)^{-2L(0)}v,
x/(1-z_{0}x))\nonumber\\
&=&Y(e^{-z_{0}(1+z_{0}x)L(1)}(1+z_{0}x)^{-2L(0)}
e^{z_{0}(1-z_{0}x)L(1)}(1-z_{0}x)^{-2L(0)}v,x).\nonumber
\end{eqnarray}
Recall (5.3.3) of [FHL]:
\begin{eqnarray}
x_{1}^{-L(0)}L(1)x_{1}^{L(0)}=x_{1}L(1).
\end{eqnarray}
{}From this we immediately get
\begin{eqnarray}\label{eL(1)L(0)conjugation}
x_{1}^{-L(0)}e^{xL(1)}x_{1}^{L(0)}=e^{xx_{1}L(1)}.
\end{eqnarray}
In view of (\ref{eL(1)L(0)conjugation}) we have
\begin{eqnarray}
e^{-z_{0}(1+z_{0}x)L(1)}(1+z_{0}x)^{-2L(0)}e^{z_{0}(1-z_{0}x)L(1)}
(1-z_{0}x)^{-2L(0)}=1,
\end{eqnarray}
hence
\begin{eqnarray}
(Y^{[z_{0}]})^{[-z_{0}]}(v,x)=Y(v,x).
\end{eqnarray}
Continuing with Proposition \ref{pdeformz} we have:

\bp{pdeformanvz}
Let $(E,Y_{E})$ be a weak $V$-module on which $L(1)$ is locally nilpotent
and let $z_{0}$ be any complex number. Then
\begin{eqnarray}\label{edeformz=old}
\Omega_{n}(E,Y_{E})=\Omega_{n}(E,Y_{E}^{[z_{0}]}).
\end{eqnarray}
Furthermore, $e^{-z_{0}L(1)}$ is an $A_{n}(V)$-isomorphism from
$\Omega_{n}(E,Y)$ to $\Omega_{n}(E,Y^{[z_{0}]})$.
\ep

\begin{proof} From (\ref{edeformzrelation}) we easily get
\begin{eqnarray}
\Omega_{n}(E,Y_{E})\subset \Omega_{n}(E,Y_{E}^{[z_{0}]}).
\end{eqnarray}
Using this and the fact that $Y_{E}=(Y_{E}^{[z_{0}]})^{[-z_{0}]}$, we get
\begin{eqnarray}
\Omega_{n}(E,Y_{E}^{[z_{0}]})\subset \Omega_{n}(E,Y_{E}).
\end{eqnarray}
This proves (\ref{edeformz=old}). The second part follows from
Proposition \ref{pdeformz} immediately.
\end{proof}

We shall use Proposition \ref{pdeformanvz} for $z_{0}=0, -1,1$.
Let $W$ and $U$ be given as before. Set
\begin{eqnarray}
E={\cal S}({\cal{D}}_{P(-1)}(W,U)).
\end{eqnarray}
In view of Lemma \ref{lomegan} we have
\begin{eqnarray}
\Omega_{n}({\cal{D}}_{P(-1)}(W,U))= \Omega_{n}(E)
\end{eqnarray}
and it follows from Corollary \ref{cnilpotentL(1)} that $L(1)$ is 
locally nilpotent on $E$, so that
we can apply Propositions \ref{pdeformz} and \ref{pdeformanvz} to $E$.

Let $W$ be a weak $V$-module. 
For homogeneous $v\in V$ and for $w\in W$, we define
\begin{eqnarray}
& &v*_{n}w\label{eanvleft}\\
&=&
\sum_{m=0}^{n}\binom{-n-1}{m}\Res_{x}x^{-n-m-1}(1+x)^{\wt v+n}Y(v,x)w,
\nonumber\\
& &
w*_{n}v\label{eanvright}\\
&=&\sum_{m=0}^{n}\binom{-n-1}{m}(-1)^{n-m}\Res_{x}x^{-n-m-1}
(1+x)^{\wt v+m-1}Y(v,x)w.\nonumber
\end{eqnarray}
Then extend the definition by linearity. 

Now, we are in a position to prove our key result.

\bp{panvbimodule}
Let $W$ be a weak $V$-module and let $U$ be a vector space. 
Let 
$$f\in \Hom (A'_{n}(W),U)=\Omega_{n}({\cal{D}}_{P(-1)}(W,U),Y^{L})
\cap \Omega_{n}({\cal{D}}_{P(-1)}(W,U),Y^{R})$$
and $w\in W$. Then 
\begin{eqnarray}
& &\left(\Res_{x}x^{\wt v}(Y^{L})^{[1]}(v,x)f\right)(w)=f(w*_{n}v)\\
& &\left(\Res_{x}x^{\wt v}(Y^{R})^{[-1]}(v,x)f\right)(w)=f(\theta(v)*_{n}w)
\end{eqnarray}
for homogeneous $v\in V$, where
\begin{eqnarray}
\theta (v)=e^{L(1)}(-1)^{L(0)}v
\end{eqnarray}
(cf. (\ref{etheta})).
\ep

\begin{proof} First,
using (\ref{eL(1)L(0)conjugation}) we get ([FHL], (5.3.1)):
\begin{eqnarray}
& &e^{xL(1)}(-x^{-2})^{L(0)}e^{x^{-1}L(1)}=(-x^{2})^{-L(0)},\\
& &e^{xL(1)}(-x^{-2})^{L(0)}e^{(x+1)^{-1}L(1)}=e^{x/(x+1)L(1)}(-x^{-2})^{L(0)}.
\end{eqnarray}
Because
\begin{eqnarray}
\left(\sum_{m=0}^{n}\binom{-n-1}{m}(-1)^{n+1-m}x^{m}\right)
(-1+x)^{n+1}\in 1+x^{n+1}{\C}[[x]],
\end{eqnarray}
for $k\ge n$,
\begin{eqnarray}
\Res_{x}x^{\wt v+k}Y^{L}(v,x)f=0.
\end{eqnarray}
Since for any homogeneous $u\in V$,
\begin{eqnarray}
(-1+x)^{\wt u+n}Y^{L}(u,x)f=(x-1)^{\wt u+n}f Y^{o}(u,x-1)
\end{eqnarray}
(Proposition \ref{pan(W)}), we have
\begin{eqnarray}
& &(-1+x)^{\wt v+n}Y^{L}(e^{(-1+x)^{-1}L(1)}v,x)f\\
&=&(x-1)^{\wt v+n}f Y^{o}(e^{(x-1)^{-1}L(1)}v,x-1),\nonumber
\end{eqnarray}
noting that $\wt L(1)^{i}v=\wt v-i$ for $i\ge 0$.
Using (\ref{edefylw0}) and all the above information we have
\begin{eqnarray}
& &\left(\Res_{x}x^{\wt v-1}(Y^{L})^{[1]}(v,x)f\right)(w)\nonumber\\
&=&\Res_{x}(-1)^{\wt v-1}x^{\wt v-1}(-1+x)^{\wt v-1}
\left(Y^{L}(e^{(-1+x)^{-1}L(1)}v,x)f\right)(w)\nonumber\\
&=&\Res_{x}\sum_{m=0}^{n}\binom{-n-1}{m}(-1)^{\wt v+n-m}x^{m+\wt v-1}
(-1+x)^{\wt v+n}\cdot \nonumber\\
& &\cdot \left(Y^{L}(e^{(-1+x)^{-1}L(1)}v,x)f\right)(w)\nonumber\\
&=&\Res_{x}\sum_{m=0}^{n}\binom{-n-1}{m}(-1)^{\wt v+n-m}x^{m+\wt v-1}
(x-1)^{\wt v+n}\cdot \nonumber\\
& &\cdot f\left(Y^{o}(e^{(x-1)^{-1}L(1)}v,x-1)w\right)
\nonumber\\
&=&\sum_{m=0}^{n}\binom{-n-1}{m}(-1)^{\wt v+n-m}\nonumber\\
& &\cdot \Res_{x}(x+1)^{m+\wt v-1}
x^{\wt v+n}f\left( Y^{o}(e^{x^{-1}L(1)}v,x)w\right)\nonumber\\
&=&\sum_{m=0}^{n}\binom{-n-1}{m}(-1)^{\wt v+n-m}\cdot \nonumber\\
& &\cdot \Res_{x}(x+1)^{m+\wt v-1}x^{\wt v+n}
f(Y(e^{xL(1)}(-x^{-2})^{L(0)}e^{x^{-1}L(1)}
v,x^{-1})w)\nonumber\\
&=&\sum_{m=0}^{n}\binom{-n-1}{m}(-1)^{\wt v+n-m}\cdot\nonumber\\
& &\cdot \Res_{x} (x+1)^{m+\wt v-1}
x^{\wt v+n}f(Y((-x^{2})^{-L(0)}v,x^{-1})w)\nonumber\\
&=&\Res_{x}\sum_{m=0}^{n}\binom{-n-1}{m}(-1)^{n-m}(x+1)^{m+\wt v-1}
x^{-\wt v+n}f(Y(v,x^{-1})w)\nonumber\\
&=&\Res_{x}\sum_{m=0}^{n}\binom{-n-1}{m}(-1)^{n-m}
(x^{-1}+1)^{m+\wt v-1}x^{\wt v-n-2}f(Y(v,x)w)\nonumber\\
&=&\Res_{x}\sum_{m=0}^{n}\binom{-n-1}{m}(-1)^{n-m}
x^{-n-m-1}(1+x)^{\wt v+m-1}f(Y(v,x)w)\nonumber\\
&=&f(w*_{n}v).\nonumber
\end{eqnarray}

Similarly, using the fact
\begin{eqnarray}
\left(\sum_{m=0}^{n}\binom{-n-1}{m}x^{m}\right)
(1+x)^{\wt v+n}\in 1+x^{n+1}{\C}[[x]]
\end{eqnarray}
we have
\begin{eqnarray*}
& &\left(\Res_{x}x^{\wt v-1}(Y^{R})^{[-1]}(v,x)f\right)(w)\nonumber\\
&=&\Res_{x}x^{\wt v-1}(1+x)^{\wt v-1}
\left(Y^{R}(e^{(1+x)^{-1}L(1)}v,x)f\right)(w)\nonumber\\
&=&\Res_{x}\sum_{m=0}^{n}\binom{-n-1}{m}x^{m+\wt v-1}
(1+x)^{\wt v+n}\left(Y^{R}(e^{(1+x)^{-1}L(1)}v,x)f\right)(w)\nonumber\\
&=&\Res_{x}\sum_{m=0}^{n}\binom{-n-1}{m}x^{m+\wt v-1}
(x+1)^{\wt v+n}f\left( Y^{o}(e^{(x+1)^{-1}L(1)}v,x)w\right)\nonumber\\
&=&\Res_{x}\sum_{m=0}^{n}\binom{-n-1}{m}x^{m+\wt v-1}
(x+1)^{\wt v+n}\cdot \nonumber\\
& &\cdot f\left(Y(e^{xL(1)}(-x^{-2})^{L(0)}e^{(x+1)^{-1}L(1)}
v,x^{-1})w\right)\nonumber\\
&=&\Res_{x}\sum_{m=0}^{n}\binom{-n-1}{m}(-1)^{\wt v}x^{m-\wt v-1}
(x+1)^{\wt v+n}\cdot \nonumber\\
& &\cdot f\left(Y(e^{x/(x+1)L(1)}(-x^{-2})^{L(0)}
v,x^{-1})w\right)\nonumber\\
&=&\Res_{x}\sum_{m=0}^{n}\binom{-n-1}{m}(-1)^{\wt v}x^{-m+\wt v-1}
(x^{-1}+1)^{\wt v+n}\cdot \nonumber\\
& &\cdot f\left(Y(e^{(1+x)^{-1}L(1)}
(-x^{2})^{L(0)}v,x)w\right)\nonumber\\
&=&\Res_{x}\sum_{m=0}^{n}\binom{-n-1}{m}(-1)^{\wt v}
x^{-n-m-1}(1+x)^{\wt v+n}
f\left(Y(e^{(1+x)^{-1}L(1)}v,x)w\right)\nonumber\\
&=&\sum_{m=0}^{n}\sum_{i\ge 0}\frac{1}{i!}
\binom{-n-1}{m}(-1)^{\wt v}\cdot\\
& &\cdot \Res_{x}x^{-n-m-1}(1+x)^{\wt (L(1)^{i}v)+n}
f\left(Y(L(1)^{i}v,x)w\right)\nonumber\\
&=&f(\theta(v)*_{n}w).
\end{eqnarray*}
This completes the proof. 
\end{proof}

One can in principle use similar arguments to those in [Z1], [FZ] and [DLM2]
to show that the left and right actions of $V$ on $W$, 
defined by (\ref{eanvleft}) and (\ref{eanvright}),
give rise to an $A_{n}(V)$-bimodule structure on $A'_{n}(W)$, or 
$A_{n}(W)$ defined below.
(From the proof of Theorem 2.3 of [DLM2], to prove the associativity 
for the right action it seems that we need to prove at least one more
combinatorial identity in addition to those proved in [DLM2].)
As a matter of fact, this easily follows from 
Proposition \ref{panvbimodule} and the  (canonical and deformed)
$A_{n}(V)\otimes A_{n}(V)$-module structures on
$$\Omega_{n}({\cal{D}}_{P(-1)}(W,U),Y^{L})
\cap \Omega_{n}({\cal{D}}_{P(-1)}(W,U),Y^{R}).$$

\bp{pan'wbimodule}
Let $W$ be a weak $V$-module. Then the left and right actions of $V$ on $W$, 
defined by (\ref{eanvleft}) and (\ref{eanvright}),
give rise to an $A_{n}(V)$-bimodule structure on $A'_{n}(W)$.
\ep

\begin{proof} Let $U={\C}$. For homogeneous $v\in V$ we set
\begin{eqnarray}
& &o^{[1]}_{L}(v)=\Res_{x}x^{\wt v-1}(Y^{L})^{[1]}(v,x),\\
& &o^{[-1]}_{R}(v)=\Res_{x}x^{\wt v-1}(Y^{R})^{[1]}(v,x).
\end{eqnarray}
Then extend the definition by linearity.
It follows from Theorem \ref{tDWU} and 
Propositions \ref{pomeganw}, \ref{pdeformz} and \ref{pdeformanvz}
that $o^{[1]}_{L}\otimes o^{[-1]}_{R}$ gives rise to an 
$A_{n}(V)\otimes A_{n}(V)$-structure on
$$\Omega_{n}({\cal{D}}_{P(-1)}(W,\C),Y^{L})
\cap \Omega_{n}({\cal{D}}_{P(-1)}(W,\C),Y^{R})=\Hom (A_{n}'(W),{\C}).$$
In particular,
\begin{eqnarray}\label{eleftrightann}
o^{[1]}_{L}(O_{n}(V))=o^{[-1]}_{R}(O_{n}(V))=0.
\end{eqnarray}
The following arguments are classical and routine in nature.
Let $u,v\in V$ be homogeneous and let $w\in W$.
For any $f\in \Hom (A_{n}'(W),{\C})$, because
$$o^{[1]}_{L}(v)f\in \Hom (A_{n}'(W),{\C}),$$
in view of Proposition \ref{panvbimodule} we have
\begin{eqnarray}
f((u\circ_{n}w)*_{n}v)=\<o^{[1]}_{L}(v)f,u\circ_{n}w\>=0.
\end{eqnarray}
Since $f$ is arbitrary, we must have
\begin{eqnarray}
(u\circ_{n}w)*_{n}v\in O_{n}'(W).
\end{eqnarray}
Using Proposition \ref{panvbimodule} and (\ref{eleftrightann}) we have
\begin{eqnarray}
\<f, w*_{n}(u\circ_{n}v)\>=\<o^{[1]}_{L}(u\circ_{n}v)f,w\>=0
\end{eqnarray}
for every $f\in \Hom (A_{n}'(W),{\C})$. Consequently,
\begin{eqnarray}
w*_{n}(u\circ_{n}v)\in O_{n}'(W).
\end{eqnarray}
Similarly, using the fact that $\theta$ gives rise to the
involution $\theta$ of $A_{n}(V)$ (Proposition \ref{pdlm2}) we have
\begin{eqnarray}
v*_{n}(u\circ_{n}w),\;\;\;\; (u\circ_{n}v)*_{n}w\in O_{n}'(W).
\end{eqnarray}
Then the left action and right action of $V$ on $W$ give
rise to a left action and right action of $A_{n}(V)$ on $A_{n}'(W)$.
The rest can be proved similarly.
\end{proof}

Motivated by the definition of $O_{n}(V)$ 
we define
\begin{eqnarray}
O_{n}(W)=O_{n}'(W)+(L(-1)+L(0))W.
\end{eqnarray}
The proof of of Lemma 2.1 of [DLM2] directly gives:

\bl{lzhuO(W)}
Let $W$ be a weak $V$-module, let $w\in W$ and let $v\in V$ be homogeneous. 
Then 
\begin{eqnarray}
v*_{n}w-w*_{n}v\equiv \Res_{x}(1+x)^{\wt v-1}Y(v,x)w\;\;\;\mod \; O_{n}(W).
\end{eqnarray}
\el

Set
\begin{eqnarray}
A_{n}(W)=W/O_{n}(W).
\end{eqnarray}
Then we have:

\bc{canvbimodule0}
The subspace $O_{n}(W)$ of $W$ is stable
under the left and right actions of $V$ on $W$, defined by (\ref{eanvleft}) 
and (\ref{eanvright}), and the quotient space $A_{n}(W)$ 
is an $A_{n}(V)$-bimodule which is the quotient module of $A_{n}'(W)$ modulo
$O_{n}(W)/O'_{n}(W)$.
\ec

\begin{proof} In view of Proposition \ref{pan'wbimodule} we only need to prove
\begin{eqnarray}
& &(L(-1)w+L(0)w)*_{n}v\in O_{n}(W),\\
& &v*_{n}(L(-1)w+L(0)w)\in O_{n}(W)
\end{eqnarray}
for $v\in V$, $w\in W$. Let us assume $v$ is homogeneous.
First, from the proof of Lemma 2.2 in [DLM2] we have
\begin{eqnarray}
\mbox{}\;\;\;\;\;\;\;
(L(-1)v+L(0)v)*_{n}w=(-1)^{n}(2n+1)\binom{2n+1}{n}(v\circ_{n} w)\in O_{n}(W).
\end{eqnarray}
Then using the fact
\begin{eqnarray}\label{efact}
[L(-1)+L(0),Y(v,x)]&=&(1+x)Y(L(-1)v,x)+Y(L(0)v,x)\\
&=&(1+x)(d/dx)Y(v,x)+(\wt v) Y(v,x),\nonumber
\end{eqnarray}
we get
\begin{eqnarray}
& &v*_{n}(L(-1)+L(0))w\\
&=&(L(-1)+L(0))(v*_{n}w)+(L(-1)v+L(0)v)*_{n}w\in O_{n}(W).\nonumber
\end{eqnarray}
Using Lemma \ref{lzhuO(W)} and (\ref{efact}),  we get
\begin{eqnarray}
& &(L(-1)w+L(0)w)*_{n}v\\
&\equiv& v*_{n}(L(-1)w+L(0)w)\nonumber\\
& &-\Res_{x}(1+x)^{\wt v-1}Y(v,x)(L(-1)w+L(0)w)
\;\;\mod\; O_{n}(W)\nonumber\\
&\equiv& -\Res_{x}(1+x)^{\wt v-1}Y(v,x)(L(-1)w+L(0)w)
\;\;\mod\; O_{n}(W)\nonumber\\
&=&-\Res_{x}(1+x)^{\wt v-1}(L(-1)+L(0))Y(v,x)w\nonumber\\
& &+\Res_{x}(1+x)^{\wt v}
(d/dx)Y(v,x)w+\Res_{x}(\wt v) (1+x)^{\wt v-1}Y(v,x)w\nonumber\\
&=&-\Res_{x}(1+x)^{\wt v-1}(L(-1)+L(0))Y(v,x)w\nonumber\\
&\equiv&0\;\;\mod\; O_{n}(W).\nonumber
\end{eqnarray}
(This argument is also similar to one in the proof Lemma 2.2 of [DLM2].)
\end{proof}

Because $A_{n}'(W)$ is an $A_{n}(V)$-bimodule 
and $\theta$ is 
an involution of $A_{n}(V)$, from the classical fact
$\Hom (A_{n}'(W),U)$ becomes an $A_{n}(V)\otimes A_{n}(V)$-module with
\begin{eqnarray}\label{etensormodule}
((a_{1},a_{2})f)(w)=f(\theta(a_{2})wa_{1})
\end{eqnarray}
for $a_{1}, a_{2}\in A_{n}(V),\; f\in \Hom (A_{n}'(W),U),\; w\in A_{n}(V)$. 
We refer to this $A_{n}(V)\otimes A_{n}(V)$-module structure as 
the {\em canonical dual module structure}.

Combining Propositions \ref{panvbimodule} with \ref{pdeformanvz}
we immediately have:

\bt{trelation}
Let $W$ be a weak $V$-module and $U$ a vector space. Let 
$$\eta: \Hom (A_{n}'(W),U)\rightarrow 
\Omega_{n}({\cal{D}}_{P(-1)}(W,U),Y^{L})\cap 
\Omega_{n}({\cal{D}}_{P(-1)}(W,U),Y^{R})$$
be the natural identification map (Proposition \ref{pan(W)}).
Then the linear map 
$$\sigma:=e^{L^{R}(1)-L^{L}(1)}\circ \eta$$
is an $A_{n}(V)\otimes A_{n}(V)$-isomorphism where $\Hom (A_{n}'(W),U)$
is equipped with the canonical dual $A_{n}(V)\otimes A_{n}(V)$-module 
structure and 
$$\Omega_{n}({\cal{D}}_{P(-1)}(W,U),Y^{L})\cap 
\Omega_{n}({\cal{D}}_{P(-1)}(W,U),Y^{R})$$
is equipped with the canonical module structure.
\et

Let $U$ be an $A_{n}(V)$-module. Then $U=\Hom_{A_{n}(V)}(A_{n}(V),U)$.
We also have the following $A_{n}$-module inclusion relations:
$$\Hom_{A_{n}(V)}(A_{n}(V),U)\subset \Hom_{A_{n}(V)}(A'_{n}(V),U)
\subset \Hom _{\C}(A'_{n}(V),U).$$
With Theorem \ref{trelation} we may and we should identify
$U$ as a submodule of the $A_{n}(V)$-module
$\Omega_{n}({\cal{D}}_{P(-1)}(V,U), Y^{L})$.

\bd{dinducedanvmodule}
{\em Let $U$ be an $A_{n}(V)$-module. We define $\Ind_{A_{n}(V)}^{V}U$ to be
the submodule of $({\cal{D}}_{P(-1)}(V,U),Y^{L})$, generated by
$U$ $(=\Hom_{A_{n}(V)}(A_{n}(V),U))$.}
\ed

Using the proof of Lemma 3.14 in [Li3] with some minor changes, we 
have:

\bp{pgeneratinganv}
Let $W$ be a weak $V$-module and let $U$ be an irreducible 
$A_{n}(V)$-submodule of $\Omega_{n}(W)$.  
Then the weak submodule $M$ of $W$ generated by $U$
is a lowest weight generalized $V$-module 
such that $M_{(h)}=U$ for some $h\in {\C}$ and $M_{(k+h)}=0$ for $k<-n$.
In particular, if $U$ is an irreducible 
$A_{n}(V)$-module, then $\Ind_{A_{n}(V)}^{V}\; U$ 
is a lowest weight generalized $V$-module
with $U$ being the homogeneous subspace of some weight $h$ such that
the lowest weight is no smaller than $h-n$. 
\ep

\end{document}